# BLOCKED REGULAR FRACTIONAL FACTORIAL DESIGNS WITH MINIMUM ABERRATION[1]

By Hongquan Xu

*University of California, Los Angeles*

This paper considers the construction of minimum aberration (MA) blocked factorial designs. Based on coding theory, the concept of minimum moment aberration due to Xu [*Statist. Sinica* **13** (2003) 691–708] for unblocked designs is extended to blocked designs. The coding theory approach studies designs in a row-wise fashion and therefore links blocked designs with nonregular and supersaturated designs. A lower bound on blocked wordlength pattern is established. It is shown that a blocked design has MA if it originates from an unblocked MA design and achieves the lower bound. It is also shown that a regular design can be partitioned into maximal blocks if and only if it contains a row without zeros. Sufficient conditions are given for constructing MA blocked designs from unblocked MA designs. The theory is then applied to construct MA blocked designs for all 32 runs, 64 runs up to 32 factors, and all 81 runs with respect to four combined wordlength patterns.

**1. Introduction.** Fractional factorial designs are widely used in scientific and industrial experiments. Blocking is an effective method for reducing systematic variations and therefore increasing precision of effect estimation. Experimenters often face the practical problem of choosing good fractional factorial designs and blocking schemes.

Fractional factorial designs are typically chosen according to the *minimum aberration* (MA) criterion [12], which includes the *maximum resolution* criterion [1] as a special case. The study of blocking in fractional factorial designs is complicated by the presence of two defining contrast subgroups, one for defining the fraction and another for defining the blocking scheme, therefore, resulting in two types of wordlength patterns, one for treatment

Received May 2005; revised November 2005.
[1]Supported in part by NSF Grants DMS-02-04009 and DMS-05-05728.

and another for block. The MA criterion can be applied to the treatment and block wordlength patterns separately. However, MA designs with respect to one wordlength pattern may not have MA with respect to the other wordlength pattern. One approach, that taken by Sun, Wu and Chen [19] and Mukerjee and Wu [15], is to consider the concept of admissible blocking schemes, but it often leads to too many admissible designs. Another approach is to combine the treatment and block wordlength patterns into one single wordlength pattern so that the criterion of MA can be applied to it in the usual way. Sitter, Chen and Feder [17], Chen and Cheng [5] and Cheng and Wu [10] have proposed four combined sequences, resulting in four MA criteria (to be defined later). See [9] for a related approach.

A practical and important issue is how to construct MA blocked designs with respect to one or more criteria. This question is not adequately addressed in the literature. Most of the existing MA blocked designs rely on the work of Sun, Wu and Chen [19], who obtained the complete catalog of blocked designs with 8, 16, 32, 64 and 128 runs for up to nine factors. MA criteria rank blocked designs according to the treatment and block wordlength patterns, which are often obtained by counting words in the treatment defining contrast subgroups and alias sets. When the number of factors is large, there are a huge number of words to be counted, causing considerable difficulties in computation. For example, when a design with 64 runs and 25 factors is arranged in 8 blocks, there are $2^{22} - 1 = 4{,}194{,}303$ words to be counted. It is cumbersome and sometimes even impossible to do so for thousands or millions of different designs. This calls for alternative computational methods.

To avoid the aforementioned computational difficulties, we take a coding theory approach and propose new methods to compare and rank blocked designs without using defining contrast subgroups and alias sets. The idea is originally due to Xu [22], who proposed the concept of *minimum moment aberration* and established its equivalence to MA for unblocked designs. We extend the concept of minimum moment aberration to blocked designs for three of the four MA criteria in Section 2.

To further ease the computation burden, we study relationships among MA blocked designs under different criteria and develop a general theory on MA blocked designs. The coding theory approach studies designs in a row-wise fashion and therefore links blocked designs with nonregular and supersaturated designs. Results on nonregular and supersaturated designs are used to establish an important lower bound on blocked wordlength pattern. It is shown that a blocked design has MA with respect to all four criteria if it originates from an unblocked MA design and achieves the lower bound. It is also shown that a regular design can be partitioned into maximal blocks if and only if it contains a row (i.e., treatment combination) without zeros. Sufficient conditions are given for constructing MA blocked designs from unblocked MA designs. Some technical lemmas are presented in Section 3 and the main results are given in Section 4. We shall point out that, for simplicity, we focus entirely on regular designs, even though most of the results can be easily extended to nonregular designs.

With the concept of minimum moment aberration and developed theory, we present methods to construct MA blocked designs in Section 5. We obtain MA blocked designs for all 8, 16, 27 and 32 runs, 64 runs up to 32 factors, and all 81 runs with respect to four combined wordlength patterns. The difference among MA blocked designs under different criteria is summarized.

The rest of this section introduces some background. A *regular* $s^{n-k}$ *design* is defined by $k$ treatment defining words, which form the *treatment defining contrast subgroup*. The *resolution* [1] is the length of the shortest word in the treatment defining contrast subgroup. For $i = 1, \ldots, n$, let $A_{i,0}$ denote the number of words of length $i$ in its treatment defining contrast subgroup. For two unblocked regular $s^{n-k}$ designs $D_1$ and $D_2$, let $r$ be the smallest integer such that $A_{r,0}(D_1) \neq A_{r,0}(D_2)$. Then $D_1$ is said to have less aberration than $D_2$ if $A_{r,0}(D_1) < A_{r,0}(D_2)$. If there is no design with less aberration than $D_1$, then $D_1$ has MA. In short, the MA criterion sequentially minimizes $A_{1,0}, A_{2,0}, \ldots, A_{n,0}$.

To arrange a regular $s^{n-k}$ design in $s^p$ blocks of size $s^{n-k-p}$, one can choose $p$ independent block defining words, which form the *block defining contrast subgroup*. There are $(s^p - 1)/(s - 1)$ block effects, each confounded with $s^k$ treatment effects. For $i = 1, \ldots, n$, let $A_{i,1}$ denote the number of treatment words of length $i$ that are confounded with some block effects.

As done in the literature, we shall only consider regular main effect (RME) designs where none of the main effects is aliased with another main effect or confounded with a block effect. It is evident that, for RME designs, $A_{1,0} = A_{2,0} = A_{1,1} = 0$. The vectors $W_t = (A_{3,0}, \ldots, A_{n,0})$ and $W_b = (A_{2,1}, \ldots, A_{n,1})$ are called the treatment and block wordlength pattern, respectively. Let $A_{0,1} = 0$ for convenience.

MA criteria for blocked designs differ in how the treatment and block wordlength patterns are combined. Sitter, Chen and Feder [17] first proposed the combined wordlength pattern

$$(1) \qquad W_{\text{scf}} = (A_{3,0}, A_{2,1}, A_{4,0}, A_{3,1}, A_{5,0}, A_{4,1}, \ldots),$$

where $A_{i,1}$ is ranked after $A_{i+1,0}$ for $i = 2, 3, \ldots$. Chen and Cheng [5] pointed out that the ordering of wordlength patterns in (1) violates the hierarchical assumption, and proposed the sequence

$$(2) \qquad W_{\text{cc}} = (3A_{3,0} + A_{2,1}, A_{4,0}, 10A_{5,0} + A_{3,1}, A_{6,0}, \ldots),$$

where the sum of $\binom{2i-1}{i} A_{2i-1,0}$ and $A_{i,1}$ is ranked before $A_{2i,0}$ for $i = 2, 3, \ldots$. Cheng and Wu [10] proposed the two combined wordlength patterns

$$(3) \qquad W_1 = (A_{3,0}, A_{4,0}, A_{2,1}, A_{5,0}, A_{6,0}, A_{3,1}, \ldots),$$

(4) $$W_2 = (A_{3,0}, A_{2,1}, A_{4,0}, A_{5,0}, A_{3,1}, A_{6,0}, \ldots),$$

where $A_{i,1}$ is ranked after $A_{2i,0}$ in $W_1$ and after $A_{2i-1,0}$ in $W_2$ for $i = 2, 3, \ldots$. We shall mention that sequence (4) was first proposed by Chen and Cheng [5] and later independently by Zhang and Park [25] and Cheng and Wu [10].

Four MA criteria result from sequentially minimizing the corresponding combined wordlength patterns. MA blocked designs under the $W$ sequence are called MA $W$ designs.

An *orthogonal array* (OA) of $N$ runs, $n$ columns, $s$ levels and strength $t$, denoted by $OA(N, n, s, t)$, is an $N \times n$ matrix in which all possible $s^t$ level combinations appear equally often as rows for any set of $t$ columns.

**2. A coding theory approach: minimum moment aberration.** For a prime power $s$, let $GF(s)$ be the finite field of $s$ elements. Let $V_n$ be the $n$-dimensional row vector space over $GF(s)$, that is, $V_n = \{(v_1, \ldots, v_n) : v_i \in GF(s) \text{ for } i = 1, \ldots, n\}$.

An $[n, m]$ *linear code* over $GF(s)$ is a vector subspace of $V_n$ with dimension $m$ so that it has $s^m$ distinct vectors. An $[n, m]$ linear code $D$ can be specified by an $m \times n$ *generator matrix* $G$ whose rows form a basis for the code. Then $D = \{u \in V_n : u = vG, v \in V_m\}$. A regular $s^{n-k}$ design is an $[n, n-k]$ linear code over $GF(s)$. For an introduction to coding theory, see [13], Chapter 4, and [20].

Consider arranging a regular $s^{n-k}$ design in $s^p$ equal-sized blocks. A design of this kind is called a regular $(s^{n-k} : s^p)$ design. Such a design is specified by a pair of matrices $T$ and $B$, defined over $GF(s)$ and of orders $(n-k) \times n$ and $(n-k) \times p$, respectively, such that $T$ has full row rank and $B$ has full column rank. Then a typical block of the design consists of all level combinations of the form $uT$, with $u \in V_{n-k}$ and $uB = v$, where $v$ is any fixed vector in $V_p$. Different blocks correspond to different choices of $v$. Since $B$ has full column rank $p$, there are $s^p$ choices of $v$, leading to a division of the $s^{n-k}$ level combinations into $s^p$ blocks. See [9] and [15].

Let $L_p = (s^p - 1)/(s - 1)$ throughout this paper. Suppose that the columns of $B$ are $b_1, \ldots, b_p$. Let $F$ be the $(n-k) \times L_p$ matrix whose columns are $\lambda_1 b_1 + \cdots + \lambda_p b_p$, where $\lambda_i \in GF(s)$, at least one $\lambda_i \neq 0$ and the first nonzero $\lambda_i$ is 1.

The columns of $T$ and $F$ can be viewed as points of $PG(n-k-1, s)$, the projective geometry of dimension $n-k-1$ over $GF(s)$. In the terminology of projective geometry, $F$ is a $(p-1)$-flat in $PG(n-k-1, s)$. Then a regular $(s^{n-k} : s^p)$ design is an RME design if and only if $T$ and $F$ are disjoint; see [5] and [15].

Let $G = (T, F)$ be the $(n-k) \times (n + L_p)$ matrix and $D$ be the linear code generated by $G$. For convenience, write $D = (D_T, D_F)$, where $D_T$ is the $N \times n$ treatment matrix and $D_F$ is the $N \times L_p$ block matrix, with $N = s^{n-k}$. For integers $t \geq 0$, define moments

$$K_{t,0}(D) = N^{-2} \sum_{i=1}^{N} \sum_{j=1}^{N} [\delta_{ij}(D_T)]^t, \tag{5}$$

$$K_{t,1}(D) = N^{-2} \sum_{i=1}^{N} \sum_{j=1}^{N} [\delta_{ij}(D_T)]^t \delta_{ij}(D_F), \tag{6}$$

where $\delta_{ij}(D_T)$ and $\delta_{ij}(D_F)$ are the number of coincidences between the $i$th and $j$th rows of $D_T$ and $D_F$, respectively. For two vectors $u = (u_1, \ldots, u_n)$ and $v = (v_1, \ldots, v_n)$, the number of coincidences is the number of $i$'s such that $u_i = v_i$. We take $0^0 = 1$ throughout the paper.

REMARK 1. The definitions of $K_{t,0}(D)$ and $K_{t,1}(D)$ given in (5) and (6) work for both regular and nonregular designs. For regular designs, the double summation can be replaced with a single summation; for example, (6) can be simplified to

$$K_{t,1}(D) = N^{-1} \sum_{i=1}^{N} [\delta_{ij}(D_T)]^t \delta_{ij}(D_F), \tag{7}$$

where $j$ can be any row number.

REMARK 2. Note that $D_F$ is a replicated $OA(s^p, L_p, s, 2)$. It follows from Lemma 1 of [14] that $\delta_{ij}(D_F)$ takes on only two different values. Specially, let $y_1, \ldots, y_N$ be the rows of $D_F$. Then

$$\delta_{ij}(D_F) = \begin{cases} L_p = (s^p - 1)/(s - 1), & \text{if } y_i = y_j, \\ L_{p-1} = L_p - s^{p-1}, & \text{otherwise.} \end{cases} \tag{8}$$

For an integer $k$, let $\binom{x}{k} = x(x-1)\cdots(x-k+1)/k!$ if $k > 0$, $\binom{x}{0} = 1$ and $\binom{x}{k} = 0$ if $k < 0$. For integers $k, j \geq 0$, let $S(k,j)$ be a Stirling number of the second kind, that is, the number of ways of partitioning a set of $k$ elements into $j$ nonempty sets. It is well known that $S(k,j) = (1/j!) \sum_{i=0}^{j} (-1)^{j-i} \binom{j}{i} i^k$ for $k \geq j \geq 0$. For integers $k, i \geq 0$, define

$$Q_k(i; n, s) = (-1)^i \sum_{j=0}^{k} j! S(k,j) s^{-j} (s-1)^{j-i} \binom{n-i}{j-i}. \tag{9}$$

For integers $t, i \geq 0$, define

$$c_t(i; n, s) = (s-1) \sum_{k=0}^{t} (-1)^k \binom{t}{k} n^{t-k} Q_k(i; n, s). \tag{10}$$

It is easy to show that $S(k,k) = 1$, $Q_k(k;n,s) = (-1)^k s^{-k} k!$ and $Q_k(i;n,s) = 0$ when $i > k$. Therefore, $c_t(t;n,s) = s^{-t}(s-1)t!$ and $c_t(i;n,s) = 0$ when $i > t$.

The following two lemmas regarding unblocked designs are from Xu [22, 23].

LEMMA 1. *For a regular $s^{n-k}$ design $D$ and integers $t \geq 0$,*

$$K_{t,0}(D) = \sum_{i=0}^{\min(t,n)} c_t(i;n,s) A_{i,0}(D), \tag{11}$$

*where $c_t(i;n,s)$ are constants defined in* (10) *and $A_{0,0}(D) = 1/(s-1)$.*

LEMMA 2. *Sequentially minimizing $K_{1,0}, K_{2,0}, \ldots, K_{n,0}$ is equivalent to sequentially minimizing $A_{1,0}, A_{2,0}, \ldots, A_{n,0}$.*

The *minimum moment aberration* criterion [22] sequentially minimizes $K_{1,0}, K_{2,0}, \ldots, K_{n,0}$. Lemma 2 implies that the minimum moment aberration criterion is equivalent to the MA criterion for unblocked designs.

Extending Lemma 1 to blocked designs, we have the following result.

THEOREM 1. *For a regular $(s^{n-k} : s^p)$ design $D$ and integers $t \geq 0$,*

$$K_{t,1}(D) = s^{-1} \sum_{i=0}^{\min(t,n)} c_t(i;n,s)[A_{i,1}(D) + L_p A_{i,0}(D)], \tag{12}$$

*where $L_p = (s^p - 1)/(s - 1)$, $c_t(i;n,s)$ are constants defined in* (10) *and $A_{0,0}(D) = 1/(s-1)$.*

The proof of Theorem 1 requires the generalized *Pless power moment identities*, a fundamental result in coding theory. For clarity, all proofs are given in the Appendix.

For an RME design $D$, $A_{1,0}(D) = A_{2,0}(D) = A_{0,1}(D) = A_{1,1}(D) = 0$. From (11) and (12), we obtain $K_{1,0}(D) = s^{-1}n$, $K_{2,0}(D) = s^{-2}n(n+s-1)$, $K_{0,1}(D) = s^{-1}L_p$ and $K_{1,1}(D) = s^{-2}nL_p$. Furthermore,

$$K_{3,0}(D) = 6s^{-3}(s-1)A_{3,0}(D) + s^{-3}n(n^2 + 3ns + s^2 - 3n - 3s + 2), \tag{13}$$

$$K_{2,1}(D) = 2s^{-3}(s-1)A_{2,1}(D) + s^{-3}n(n+s-1)L_p. \tag{14}$$

We can define three minimum moment aberration criteria for blocked designs by replacing $A_{i,0}$ and $A_{i,1}$ with $K_{i,0}$ and $K_{i,1}$ in (1), (3) and (4). Because $c_t(t;n,s)$ is a positive constant, it follows from (11) and (12) that the minimum moment aberration criterion with respect to $W_{\text{scf}}$, $W_1$ or $W_2$ is equivalent to its corresponding MA criterion.

The MA $W_{cc}$ criterion defined in (2) is more complicated than the other three criteria. Nevertheless, from (13) and (14), we obtain

$$K_{3,0}(D) + K_{2,1}(D) = 2s^{-3}(s-1)[3A_{3,0}(D) + A_{2,1}(D)]$$
$$+ s^{-3}n[(n^2 + 3ns + s^2 - 3n - 3s + 2) + (n+s-1)L_p].$$

Therefore, minimizing $K_{3,0} + K_{2,1}$ is equivalent to minimizing $3A_{3,0} + A_{2,1}$.

**3. Some lemmas.** Suppose that $D = (D_T, D_F)$ is a regular $(s^{n-k} : s^p)$ design. Let $x_1, \ldots, x_N$ be the rows of $D_T$ and $y_1, \ldots, y_N$ be the rows of $D_F$, where $N = s^{n-k}$. For $m = 1, \ldots, s^p$, let $D_m$ be the $s^{n-k-p} \times n$ treatment matrix corresponding to the $m$th block. For integers $t \geq 0$, define moments

$$K_t(D_m) = s^{-2(n-k-p)} \sum_{i=1}^{s^{n-k-p}} \sum_{j=1}^{s^{n-k-p}} [\delta_{ij}(D_m)]^t,$$

where $\delta_{ij}(D_m)$ is the number of coincidences between the $i$th and $j$th rows of $D_m$. Let $B_m = \{i : x_i \text{ is a row of } D_m, 1 \leq i \leq N\}$. It is evident that $i \in B_m$ and $j \in B_m$ for some $m$ if and only if $y_i = y_j$. It is useful to express $K_t(D_m)$ in terms of the original design $D_T$ as

(15) $$K_t(D_m) = s^{-2(n-k-p)} \sum_{i \in B_m} \sum_{j \in B_m} [\delta_{ij}(D_T)]^t.$$

Without loss of generality, assume that $D_1$ contains the null treatment (i.e., a row of zeros) and call $D_1$ the *principal block*. Then $D_1$ is an $[n, n-k-p]$ linear code over $GF(s)$ and other blocks $D_m$, $2 \leq m \leq s^p$, are cosets of $D_1$; therefore,

(16) $$K_t(D_m) = K_t(D_1) \quad \text{for } m = 2, \ldots, s^p.$$

Note that $D_1$ is possibly a supersaturated design in which the number of columns is larger than the number of rows.

The next result shows that $K_{t,1}(D)$ is determined by $K_{t,0}(D)$ and $K_t(D_1)$.

LEMMA 3. *Suppose that $D$ is a regular $(s^{n-k} : s^p)$ design and $D_1$ is its principal block. For integers $t \geq 0$, $K_{t,1}(D) = L_{p-1}K_{t,0}(D) + s^{-1}K_t(D_1)$.*

LEMMA 4. *Suppose that $D$ is an $(s^{n-k} : s^p)$ RME design and $D_1$ is its principal block. Let $J = n(s^{n-k-p-1} - 1)(s^{n-k-p} - 1)^{-1}$ and $\eta$ be the fractional part of $J$.*

*(i) $K_1(D_1) = s^{-1}n$ and $D_1$ is an $OA(s^{n-k-p}, n, s, 1)$.*
*(ii) $K_2(D_1) \geq s^{-(n-k-p)}[n^2 + (s^{n-k-p} - 1)(J^2 + \eta(1 - \eta))]$. The equality holds if and only if the difference among all $\delta_{ij}(D_1)$, $i < j$, does not exceed one.*

(iii) $K_{2,1}(D) \geq L_{p-1}s^{-2}n(n+s-1)+s^{-(n-k-p+1)}[n^2+(s^{n-k-p}-1)(J^2+\eta(1-\eta))]$.

A regular $(s^{n-k}:s^p)$ design $D=(D_T,D_F)$ can be viewed as an unblocked regular $s^{(n+L_p)-(k+L_p)}$ design. For clarity, denote this unblocked design as $D_{\mathrm{un}}$. For integers $t \geq 0$, define moments $K_t(D_{\mathrm{un}}) = N^{-2}\sum_{i=1}^{N}\sum_{j=1}^{N}[\delta_{ij}(D_{\mathrm{un}})]^t$, where $N=s^{n-k}$ and $\delta_{ij}(D_{\mathrm{un}})=\delta_{ij}(D_T)+\delta_{ij}(D_F)$ is the number of coincidences between the $i$th and $j$th rows of $D_{\mathrm{un}}$. The next result shows that $K_t(D_{\mathrm{un}})$ is related to $K_{t,0}(D)$, $K_{t-1,1}(D)$, $K_{t-2,1}(D)$, $K_{t-2,0}(D)$ and so on.

LEMMA 5. *For a regular $(s^{n-k}:s^p)$ design $D$ and integers $t \geq 0$,*

$$K_t(D_{\mathrm{un}}) = K_{t,0}(D) + tK_{t-1,1}(D)$$
$$+ s^{-p+1}\sum_{r=2}^{t}\binom{t}{r}[(L_p^r - L_{p-1}^r)K_{t-r,1}(D)$$
$$- L_pL_{p-1}(L_p^{r-1}-L_{p-1}^{r-1})K_{t-r,0}(D)].$$

The following two lemmas are useful to know when MA blocked designs are the same under different criteria.

LEMMA 6. *If $D$ has MA with respect to both $W_{\mathrm{scf}}$ and $W_1$, then $D$ has MA with respect to $W_2$.*

LEMMA 7. *Suppose there exists some constant $0 \leq \alpha < 3$ such that $\alpha A_{3,0}+A_{2,1}$ is minimized for $D$. If $D$ has MA with respect to both $W_{\mathrm{scf}}$ and $W_2$, then $D$ has MA with respect to $W_{\mathrm{cc}}$.*

Lemma 7 is very useful to show the MA $W_{\mathrm{cc}}$ optimality. The condition $\alpha < 3$ is necessary; see Section 5 for counterexamples.

**4. Main results.** Lemma 4(iii) and (14) together yield a lower bound of $A_{2,1}$ as follows.

THEOREM 2. *For an $(s^{n-k}:s^p)$ RME design $D$,*

$$A_{2,1}(D) \geq [2(s-1)]^{-1}$$
$$\times \{-n(n+s-1)+s^{-(n-k-p-2)}$$
$$\times [n^2+(s^{n-k-p}-1)(J^2+\eta(1-\eta))]\},$$

*where $J=n(s^{n-k-p-1}-1)(s^{n-k-p}-1)^{-1}$ and $\eta$ is the fractional part of $J$.*

Theorem 2 plays an important role in the theoretical development and construction of MA blocked designs. The lower bound is tight for $p = n - k - 1$ and $n - k - 2$. Note that an RME design achieving the lower bound does not always have MA. When $s = 2$ and $p < n - k - 2$, the lower bound can be improved in some cases if the results in [4] and [2] are used. However, the improvement is usually negligible, noting that $A_{2,1}$ must be an integer for RME designs.

COROLLARY 1. *If $D$ has MA with respect to $W_{\text{scf}}$ and $W_2$, and $D$ achieves the lower bound in Theorem 2, then $D$ has MA with respect to $W_{\text{cc}}$.*

The next result provides a sufficient condition when MA blocked designs are the same under four criteria.

THEOREM 3. *If $D_T$ has MA among all regular $s^{n-k}$ designs and $D$ achieves the lower bound in Theorem 2, then $D$ has MA with respect to $W_{\text{scf}}$, $W_1$, $W_2$ and $W_{\text{cc}}$.*

When the lower bound in Theorem 2 is achieved, the principal block $D_1$ has minimum moment aberration among all $s^{n-k-p} \times n$ designs. Theorem 3 can be generalized as follows.

COROLLARY 2. *If $D_T$ has MA among all regular $s^{n-k}$ designs and the principal block $D_1$ has minimum moment aberration among all $s^{n-k-p} \times n$ designs, then $D$ has MA with respect to $W_{\text{scf}}$, $W_1$, $W_2$ and $W_{\text{cc}}$.*

The next result gives a simple necessary and sufficient condition when a regular design can be partitioned into maximal blocks as an RME design.

THEOREM 4. *A regular $s^{n-k}$ design containing the null treatment can be partitioned into maximal $s^{n-k-1}$ blocks as an RME design if and only if it contains a row without zeros.*

Mukerjee and Wu [15] previously studied the maximal blocking problem with a projective geometric approach. They managed to obtain a complete solution for $s^{n-1}$ and $s^{n-2}$ designs. Our approach appears to be more pleasant than theirs. Theorem 4 gives a simple answer to the question.

When $s = 2$, a row without zeros is necessarily a row of all 1's. Then a row and its fold-over forms a block. The unblocked design must be a fold-over design. A regular fold-over design is also called an *even design* [11], because it contains only words of even length. Whether or not a design is an even design can be simply checked by its wordlength pattern. The following corollary is a special case of Theorem 4.

COROLLARY 3. *A regular $2^{n-k}$ design containing the null treatment can be partitioned into maximal $2^{n-k-1}$ blocks as an RME design if and only if it is an even design.*

It is of special interest to know when an unblocked MA design can be partitioned into maximal blocks. Unblocked MA $2^{n-k}$ designs were given by Chen and Wu [8] for $k = 1, 2, 3, 4$ and by Chen [6] for $k = 5$. Combining their results and Corollary 3, we have the following result. An unblocked MA $2^{n-k}$ design can be partitioned into maximal $2^{n-k-1}$ blocks as an RME design as follows:

(i) when $k = 1$ and $n$ is even,
(ii) when $k = 2$ and $n$ is a multiple of 3,
(iii) when $k = 3$ and $n = 7t + q$ for integers $t \geq 0$ and $q = 7, 11$,
(iv) when $k = 4$ and $n = 15t + q$ for integers $t \geq 0$ and $q = 8, 12, 15, 20$,
(v) when $k = 5$ and $n = 31t + q$ for integers $t \geq 0$ and $q = 16, 21, 24, 28, 31, 37, 40, 44$.

Furthermore, it is known from coding theory that even designs are the only designs of resolution IV for $5N/16 < n \leq N/2$ with $N = 2^{n-k}$; see [3]. For such $n$, an unblocked MA $2^{n-k}$ design can always be partitioned into maximal $2^{n-k-1}$ blocks as an RME design.

To describe the next result, let $\tilde{F}$ be an $(n - k - 2)$-flat and $\tilde{T}$ be the complement of $\tilde{F}$ in $PG(n - k - 1, s)$. Let $\tilde{H}_{n-k}$ be the linear code generated by $\tilde{T}$. Note that $\tilde{H}_{n-k}$ is unique up to isomorphism. It is evident that $\tilde{H}_{n-k}$ and its projection designs (i.e., subsets of columns) can be partitioned into maximal $s^{n-k-1}$ blocks as RME designs. The reverse is also true in the following sense. If an $s^{n-k}$ design can be partitioned into $s^{n-k-1}$ blocks as an RME design, then it is isomorphic to a projection design of $\tilde{H}_{n-k}$. The next result characterizes MA $(s^{n-k} : s^{n-k-1})$ RME designs.

THEOREM 5. *If $D_T$ has MA among all projection designs of $\tilde{H}_{n-k}$, then $D_T$ can be partitioned as an $(s^{n-k} : s^{n-k-1})$ RME design $D$ that has MA with respect to $W_{\text{scf}}$, $W_1$, $W_2$ and $W_{\text{cc}}$.*

Theorem 5 shows that MA blocked $(s^{n-k} : s^{n-k-1})$ designs are the same for all four criteria when they exist. As a numeric illustration, consider $s = 3$. For 27 runs, $\tilde{H}_3$ is the unique MA $3^{9-6}$ design. According to Xu [23], for $4 \leq n < 9$, MA $3^{n-(n-3)}$ designs are projection designs of the MA $3^{9-6}$ design; therefore, they can be partitioned into maximal 9 blocks and resulting RME designs have MA with respect to all four criteria. For 81 runs, $\tilde{H}_4$ is the unique MA $3^{27-23}$ design. According to Xu [23], for $5 \leq n \leq 9$ and $12 \leq n < 27$, MA $3^{n-(n-4)}$ designs are projection designs of the MA $3^{27-23}$ design; therefore, they can be partitioned into maximal 27 blocks and resulting RME designs have MA with respect to all four criteria. For $n = 10, 11$, MA $3^{n-(n-4)}$ designs are not projection designs of the MA $3^{27-23}$ design; therefore, they cannot be partitioned as RME designs with 27 blocks. The second best designs are projection designs of the MA $3^{27-23}$ design; therefore, they can be partitioned into maximal 27 blocks and the resulting RME designs have MA with respect to all four criteria.

When $s = 2$, $\tilde{H}_{n-k}$ is an even design with resolution IV. We have the following result.

COROLLARY 4. *If a regular $2^{n-k}$ design has MA among all even designs, then it can be partitioned as a $(2^{n-k}:2^{n-k-1})$ RME design that has MA with respect to $W_{\mathrm{scf}}$, $W_1$, $W_2$ and $W_{\mathrm{cc}}$.*

Recall that a regular $(s^{n-k}:s^p)$ design $D$ can be viewed as an unblocked regular $s^{(n+L_p)-(k+L_p)}$ design $D_{\mathrm{un}}$. The next result provides a sufficient condition when an MA blocked design originates from an unblocked MA design.

THEOREM 6. *If $D_T$ has MA among all regular $s^{n-k}$ designs and the unblocked design $D_{\mathrm{un}}$ has MA among all regular $s^{(n+L_p)-(k+L_p)}$ designs, then the blocked $(s^{n-k}:s^p)$ RME design $D$ has MA with respect to $W_{\mathrm{scf}}$, $W_1$, $W_2$ and $W_{\mathrm{cc}}$.*

Theorem 6 is most useful when $p = 1$. It happens frequently that an unblocked MA $s^{n-k}$ design can be extended to an unblocked MA $s^{(n+1)-(k+1)}$ design by adding an extra column. For example, according to Chen, Sun and Wu [7] and Xu [23], for 8 and 27 runs, MA unblocked designs are in sequential order for all $n$. Whenever this happens, the extra column can be used as the block generator, and the resulting $(s^{n-k}:s^1)$ design has MA with respect to $W_{\mathrm{scf}}$, $W_1$, $W_2$ and $W_{\mathrm{cc}}$.

Theorem 6 is less useful when $p > 1$ because $D_{\mathrm{un}}$ usually does not have MA. The following result is interesting in this regard.

THEOREM 7. *If $D_T$ has MA among all regular $s^{n-k}$ designs and $D$ has MA with respect to $W_{\mathrm{scf}}$, then $D$ has MA with respect to both $W_1$ and $W_2$. If, in addition, $\alpha A_{3,0}(D) + A_{2,1}(D)$ is also minimized for some constant $0 \leq \alpha < 3$, then $D$ has MA with respect to $W_{\mathrm{cc}}$.*

Theorem 7 implies that when MA blocked designs are different under $W_{\mathrm{scf}}$, $W_1$ and $W_2$, an MA $W_{\mathrm{scf}}$ blocked design must not originate from an unblocked MA design.

**5. MA blocked designs.** MA blocked designs with respect to $W_{\text{scf}}$, $W_1$ or $W_2$ can be obtained by computing and comparing moments $K_{t,0}$ and $K_{t,1}$ for all possible blocking schemes. This is a feasible task when the number of blocking schemes is not too large; see [24] for details, where MA blocked designs for all 32 runs, 64 runs up to 32 factors, and all 81 runs with respect to $W_{\text{scf}}$, $W_1$ and $W_2$ are given.

However, this method cannot be used to construct MA $W_{\text{cc}}$ designs because there is no equivalent minimum moment aberration criterion with respect to (2). Furthermore, an essential difference exists between the MA $W_{\text{cc}}$ criterion and the other three criteria. Because the MA $W_{\text{scf}}$, $W_1$ and $W_2$ criteria minimize $A_{3,0}$ first, there is no need to search over resolution III designs whenever blocking schemes from resolution IV designs exist. However, the MA $W_{\text{cc}}$ criterion minimizes $3A_{3,0} + A_{2,1}$ first. Combining $A_{3,0}$ with $A_{2,1}$ makes it more difficult to construct MA $W_{\text{cc}}$ designs than other types of MA designs. To determine the minimum of $3A_{3,0} + A_{2,1}$, a simple strategy is to search over all resolution III designs. This requires a complete catalog of resolution III designs, but such a catalog is not available for 64-run designs.

Combining the developed theory and computer search, we obtain MA $W_{\text{cc}}$ designs for all 8, 16, 27 and 32 runs, 64 runs up to 32 factors, and all 81 runs. Previously, Chen and Cheng [5] developed a theory to characterize MA $W_{\text{cc}}$ designs in terms of their blocked residual designs and obtained MA $W_{\text{cc}}$ designs for all 8 and 16 runs and 32 runs up to 20 factors.

Here we explain how to construct MA $W_{\text{cc}}$ designs for 64 runs and $n \leq 32$ with the results of Xu and Lau [24]. First, for $p = 5$ and $6 \leq n \leq 32$, by Theorem 5, MA $(2^{n-(n-6)} : 2^5)$ designs are the same under all four criteria; therefore, MA designs given by Xu and Lau [24] have MA with respect to all four criteria. Indeed, they can be easily constructed by searching over MA projection designs of the unique even $2^{32-26}$ design. Next, for $p = 1$, because unblocked MA $2^{n-(n-6)}$ designs are in sequential order for $n = 6$–7, 8–12, 14–20 and 21–33, by Theorem 6, we obtain MA $(2^{n-(n-6)} : 2^1)$ designs with respect to all four criteria for all $6 \leq n \leq 32$ but $n = 7$, 12, 13 and 20. For $n = 7$, 12, 13 or 20, according to [24], MA $W_2$ and $W_{\text{scf}}$ designs coincide and have $A_{2,1} = 0$; therefore, by Lemma 7, MA $W_{\text{cc}}$ designs also coincide with MA $W_2$ and $W_{\text{scf}}$ designs.

The situation for $p = 2, 3, 4$ is more complicated than that for $p = 1$ and 5. We first compute the lower bounds of $A_{2,1}$ in Theorem 2, which are given in Table 1. It is evident that a lower bound can be replaced by the smallest nonnegative integer that exceeds it if it is negative or not an integer. According to Xu and Lau [24], MA $W_2$ designs achieve the modified lower bounds of $A_{2,1}$ except for the following 22 cases: $p = 2$, $n = 19$–26, 31, 32; $p = 3$, $n = 29$–32; and $p = 4$, $n = 25$–32. Furthermore, MA $W_{\text{scf}}$ and $W_2$ designs coincide except for $p = 2$ and $n = 7, 12$. Then, by Corollary 1, MA $W_{\text{scf}}$ and $W_2$ designs also have MA with respect to $W_{\text{cc}}$ except for the 24 special cases, which require additional computer search.

TABLE 1
*Lower bound of $A_{2,1}$ in Theorem 2 for 64 runs, $6 \leq n \leq 32$ and $p = 2, 3, 4$*

| $p$ | 6 | 7 | 8 | 9 | 10 | 11 | 12 | 13 | 14 | 15 | 16 | 17 | 18 | 19 |
|---|---|---|---|---|---|---|---|---|---|---|---|---|---|---|
| 2 | −1.5 | −1.5 | −1.5 | −1.5 | −1.2 | −1.3 | −0.8 | −0.7 | 0 | 0 | 1 | 1.2 | 2.3 | 2.7 |
| 3 | 0 | 0 | 1 | 1.5 | 2.5 | 3.5 | 4.5 | 6 | 7 | 9 | 10.5 | 12.5 | 14.5 | 16.5 |
| 4 | 3 | 5 | 7 | 9 | 12 | 15 | 18 | 22 | 26 | 30 | 35 | 40 | 45 | 51 |

| $p$ | 20 | 21 | 22 | 23 | 24 | 25 | 26 | 27 | 28 | 29 | 30 | 31 | 32 |
|---|---|---|---|---|---|---|---|---|---|---|---|---|---|
| 2 | 3.8 | 4.5 | 5.5 | 6.5 | 7.5 | 8.7 | 9.8 | 11.3 | 12.2 | 14 | 15 | 17 | 18.2 |
| 3 | 19 | 21 | 24 | 26.5 | 29.5 | 32.5 | 35.5 | 39 | 42 | 46 | 49.5 | 53.5 | 57.5 |
| 4 | 57 | 63 | 70 | 77 | 84 | 92 | 100 | 108 | 117 | 126 | 135 | 145 | 155 |

For the 24 special cases, Theorem 2 and Lemma 7 are again used to ease computation. Consider, for example, $p = 4$ and $n = 29$. According to Xu and Lau [24], MA $W_{\text{scf}}, W_1$ and $W_2$ designs coincide and have $A_{3,0} = 0$ and $A_{2,1} = 196$. The lower bound of $A_{2,1}$ is 126. To determine the minimum of $3A_{3,0} + A_{2,1}$, we only need search over all designs with $A_{3,0} \leq (196 - 126)/3 = 23.3$, leading to $A_{3,0} \leq 23$. This is a feasible task. A complete enumeration (to be explained later) shows that there are exactly 17 regular $2^{29-23}$ designs with $A_{3,0} \leq 23$, among which one has resolution IV. It is straightforward to verify that $2.9A_{3,0} + A_{2,1}$ has minimum 196 among all 17 $2^{29-23}$ designs with $A_{3,0} \leq 23$. Then, by Lemma 7, MA $W_{\text{cc}}$ designs coincide with MA $W_2$ and $W_{\text{scf}}$ designs.

When MA $W_2$ and $W_{\text{scf}}$ designs are different or when they do not minimize $\alpha A_{3,0} + A_{2,1}$ for all $\alpha$ with $0 \leq \alpha < 3$, Lemma 7 cannot be used; then MA $W_{\text{cc}}$ designs are determined by sequentially comparing the complete sequence in (2). Fortunately, this happens only for the following five cases: $(n, p) = (7, 2), (12, 2), (25, 4), (26, 4), (29, 3)$. For the first two cases, MA $W_{\text{cc}}$ designs coincide with MA $W_2$ designs; for the last three cases, MA $W_{\text{cc}}$ designs are different from MA $W_2$ designs, which coincide with MA $W_{\text{scf}}$ and $W_1$ designs.

Now we explain how to enumerate all $2^{29-23}$ designs with $A_{3,0} \leq 23$. Note that a 3-letter word consists of three factors and there are 29 factors in a $2^{29-23}$ design. Therefore, for any $2^{29-23}$ design with $A_{3,0} = 23$, there must exist a column appearing in at least $3 \times 23/29 = 2.4$ or 3 words of length 3. Deleting that column yields a $2^{28-22}$ design with $A_{3,0} \leq 20$. Therefore, all $2^{29-23}$ designs with $A_{3,0} \leq 23$ can be enumerated by adding a column to all $2^{28-22}$ designs with $A_{3,0} \leq 20$, which in turn can be enumerated by adding a column to all $2^{27-21}$ designs with $A_{3,0} \leq 17$. This can be done sequentially in the same way as in Chen, Sun and Wu [7] and Xu [23], as long as the number of designs is not too large at each step. We shall point out the importance of the lower bound of $A_{2,1}$ in Theorem 2. Without this bound, one has to search over all $2^{29-23}$ designs with $A_{3,0} \leq 196/3 = 65.3$. This is not a feasible task because there are more than 100,000 $2^{29-23}$ designs with $A_{3,0} \leq 65$ and it is impossible to enumerate all of them with the current method and computer.

Finally, we summarize the differences of MA blocked designs under different criteria for all 8, 16, 27, 32 runs, 64 runs up to 32 factors, and all 81 runs. We observed that MA blocked designs under all four criteria are the same in most cases. This occurs for all 8 and 27 runs, which can be easily verified with Theorems 5 and 6. When MA blocked designs under four criteria are not all the same, one of the following four situations occurs:

1. MA $W_1$, $W_2$ and $W_{cc}$ designs are the same, but they differ from MA $W_{scf}$ designs.
2. MA $W_2$, $W_{scf}$ and $W_{cc}$ designs are the same, but they differ from MA $W_1$ designs.
3. MA $W_1$, $W_2$ and $W_{scf}$ designs are the same, but they differ from MA $W_{cc}$ designs.
4. MA $W_2$ and $W_{cc}$ designs are the same, but they differ from MA $W_1$ or $W_{scf}$ designs.

Situation 1 occurs once for 32 runs with $(n,p) = (6,1)$ and once for 64 runs with $(n,p) = (7,2)$, and does not occur for 16 and 81 runs. Situation 2 occurs twice for 16 runs with $(n,p) = (5,1),(5,2)$, 12 times for 32 runs, 35 times for 64 runs, and twice for 81 runs with $(n,p) = (11,1),(11,2)$. Situation 3 occurs once for 32 runs with $(n,p) = (13,3)$, three times for 64 runs with $(n,p) = (25,4),(26,4),(29,3)$, and three times for 81 runs with $(n,p) = (9,2),(17,2),(21,2)$. Situation 4 occurs only once for 64 runs with $(n,p) = (12,2)$.

Except for situation 3, MA $W_{cc}$ designs coincide with MA $W_2$ designs, which are given by Xu and Lau [24]. Table 2 gives MA designs for situation 3 with treatment and block columns in the same fashion as Cheng and Wu [10] and Xu and Lau [24]. The designs are labeled as $n$-$k \cdot i/\mathrm{B}p(W)$, where $i$ denotes the rank of the unblocked $s^{n-k}$ design under the MA criterion, $p$ denotes the number of block variables, and $W$ denotes the MA $W$-criterion. See Xu and Lau [24] for generator matrices and column labels. To save space, in Table 2 independent columns are omitted in the treatment columns; only generators are given in the block columns, treatment wordlength pattern is truncated as $W_t = (A_{3,0}, A_{4,0}, A_{5,0}, A_{6,0})$ and block wordlength pattern is truncated as $W_b = (A_{2,1}, A_{3,1}, A_{4,1}, A_{5,1})$. The last two columns in Table 2 give the numbers of clear main effects ($C1$) and of clear two-factor interactions ($C2$). A main effect or two-factor interaction is *clear* if it is not aliased with any other main effect or two-factor interaction and is not confounded with any block effect [19].

TABLE 2
*MA blocked designs for situation* 3

| Design | Treatment | $W_t$ | Block | $W_b$ | $C1$ | $C2$ |
|---|---|---|---|---|---|---|
| | | 32 runs | | | | |
| 13-8.1/B3(*) | 31 7 11 21 25 13 14 19 | 0 55 0 96 | 3 5 17 | 36 0 310 0 | 13 | 0 |
| 13-8.4/B3($W_{cc}$) | 31 7 11 21 13 14 26 3 | 4 39 32 48 | 5 10 19 | 22 76 124 288 | 4 | 0 |
| | | 64 runs | | | | |
| 25-19.1/B4(*) | 31 35 13 52 14 55 37 61 11 19 21 44 7 62 25 49 22 41 38 | 0 435 0 5440 | 3 5 9 48 | 144 0 5923 0 | 25 | 0 |
| 25-19.17/B4($W_{cc}$) | 31 35 13 52 14 55 21 37 11 19 25 38 7 26 49 22 28 50 9 | 8 378 336 4032 | 3 5 17 41 | 92 568 2688 13104 | 8 | 0 |
| 26-20.1/B4(*) | 31 35 13 52 14 55 37 61 11 19 21 44 7 62 25 49 22 41 38 26 | 0 515 0 7062 | 3 5 9 48 | 156 0 6999 0 | 26 | 0 |
| 26-20.50/B4($W_{cc}$) | 31 35 13 52 14 55 21 37 11 19 25 38 7 26 49 22 28 50 9 33 | 16 386 672 4368 | 3 5 17 41 | 100 632 3248 15960 | 0 | 0 |
| 29-23.1/B3(*) | 31 35 13 52 14 55 37 61 11 19 21 44 7 62 25 49 22 41 38 26 28 42 47 | 0 819 0 14560 | 5 17 33 | 91 0 5187 0 | 29 | 0 |
| 29-23.4/B3($W_{cc}$) | 31 35 13 52 14 55 37 61 11 19 21 44 7 62 25 49 22 41 26 28 42 56 3 | 12 707 640 11536 | 9 20 38 | 46 484 2252 14016 | 4 | 0 |
| | | 81 runs | | | | |
| 9-5.1/B2(*) | 22 9 24 31 34 | 0 18 36 12 | 4 20 | 9 30 117 162 | 9 | 0 |
| 9-5.2/B2($W_{cc}$) | 22 9 24 31 3 | 1 18 27 28 | 6 18 | 6 44 90 186 | 6 | 5 |
| 17-13.1/B2(*) | 22 9 24 31 3 25 13 37 6 18 7 35 12 | 20 336 1014 4 15 5072 | | 40 210 2079 9256 | 0 | 0 |
| 17-13.2/B2($W_{cc}$) | 22 9 24 31 3 25 13 37 6 18 7 35 16 | 23 306 1107 12 15 4952 | | 28 303 1782 9814 | 0 | 0 |
| 21-17.1/B2(*) | 22 9 24 31 3 25 13 37 6 18 7 35 12 38 15 16 19 | 51 729 3717 4 26 21819 | | 48 550 4590 32418 | 0 | 0 |
| 21-17.2/B2($W_{cc}$) | 22 9 24 31 3 25 13 37 15 23 16 34 6 38 7 18 26 | 52 720 3735 11 30 21876 | | 45 573 4545 32310 | 0 | 0 |

NOTES. (*) MA $W_{scf}$, $W_1$ and $W_2$ designs.

For all designs given in Table 2, MA $W_{cc}$ designs have a larger $A_{3,0}$ value but a smaller $A_{2,1}$ value than corresponding MA designs under the other three criteria. Indeed, these MA $W_{cc}$ designs achieve the lower bound of $A_{2,1}$ in Theorem 2, whereas MA designs under the other three criteria originate from unblocked MA designs.

Note that for 81 runs, MA $W_{scf}, W_1$ and $W_2$ designs 9-5.1/B2 have the same $3A_{3,0} + A_{2,1}$ value as the MA $W_{cc}$ design 9-5.2/B2. This also happens with 21-17.1/B2 and 21-17.2/B2. These examples show that the condition $\alpha < 3$ in Lemma 7 is necessary.

## APPENDIX

Further notation and results in coding theory are necessary in order to prove Theorem 1. The *Hamming weight* of a vector $u = (u_1, \ldots, u_n)$, denoted by $wt(u)$, is the number of its nonzero components $u_i$.

Associated with any linear code $D$ is another linear code, called its *dual* and denoted by $D^\perp$. Suppose $D$ is an $[n, m]$ linear code with generator matrix $G$ over $GF(s)$. The dual $D^\perp$ is the null space of $G$, that is, $D^\perp = \{u \in V_n : uG' = 0\}$, where $G'$ is the transpose of $G$. The dual $D^\perp$ is indeed the defining contrast subgroup of $D$.

Suppose $D$ is an $[n_1 + n_2, m]$ linear code over $GF(s)$ and $D^\perp$ is its dual code. Each vector $u$ in $D$ and $D^\perp$ can be written as $u = (u_1, u_2)$, where $u_1 \in V_{n_1}$ and $u_2 \in V_{n_2}$. Let $B_{i_1,i_2}(D)$ and $B_{i_1,i_2}(D^\perp)$ be the number of vectors in $D$ and respectively in $D^\perp$ with $wt(u_1) = i_1$ and $wt(u_2) = i_2$.

The following result, a special case of Lemma 4.3 of Xu [21], generalizes the Pless power moment identities [16].

LEMMA A.1.   *For integers* $k_1, k_2 \geq 0$,

$$s^{-m} \sum_{i_1=0}^{n_1} \sum_{i_2=0}^{n_2} i_1^{k_1} i_2^{k_2} B_{i_1,i_2}(D) = \sum_{j_1=0}^{n_1} \sum_{j_2=0}^{n_2} B_{j_1,j_2}(D^\perp) Q_{k_1}(j_1; n_1, s) Q_{k_2}(j_2; n_2, s),$$

*where* $Q_k(j; n, s)$ *is defined in* (9).

PROOF OF THEOREM 1.   Let $N = s^{n-k}$ and $n_2 = L_p$. Then $D = (D_T, D_F)$ is an $[n + n_2, n - k]$ linear code. Let $D^\perp$ be the dual code of $D$. Each vector $u$ in $D$ and $D^\perp$ can be written as $u = (u_1, u_2)$, where $u_1 \in V_n$ and $u_2 \in V_{n_2}$. It is known that the wordlength patterns are proportional to the split weight distributions of $D^\perp$ as follows: for $i = 0, \ldots, n$,

(A.1)   $A_{i,0}(D) = B_{i,0}(D^\perp)/(s-1)$   and   $A_{i,1}(D) = B_{i,1}(D^\perp)/(s-1);$

see [18] and [5]. By (7),

$$\begin{aligned} K_{t,1}(D) &= N^{-1} \sum_{i_1=0}^{n} \sum_{i_2=0}^{n_2} (n - i_1)^t (n_2 - i_2) B_{i_1,i_2}(D) \\ &= N^{-1} \sum_{i_1=0}^{n} \sum_{i_2=0}^{n_2} \sum_{k=0}^{t} \binom{t}{k} (-1)^k n^{t-k} i_1^k (n_2 - i_2) B_{i_1,i_2}(D) \\ &= N^{-1} \sum_{k=0}^{t} \binom{t}{k} (-1)^k n^{t-k} \sum_{i_1=0}^{n} \sum_{i_2=0}^{n_2} (i_1^k n_2 - i_1^k i_2) B_{i_1,i_2}(D). \end{aligned}$$

By Lemma A.1,

$$K_{t,1}(D) = \sum_{k=0}^{t} \binom{t}{k}(-1)^k n^{t-k} \sum_{j_1=0}^{n}\sum_{j_2=0}^{n_2} B_{j_1,j_2}(D^\perp)Q_k(j_1;n,s)$$
$$\times [Q_0(j_2;n_2,s)n_2 - Q_1(j_2;n_2,s)].$$

Recall that $Q_k(j;n,s) = 0$ for $j > k$. Then

$$K_{t,1}(D) = \sum_{k=0}^{t}\binom{t}{k}(-1)^k n^{t-k}\sum_{j_1=0}^{n} Q_k(j_1;n,s)\Delta(D^\perp,j_1;n_2,s),$$

where $\Delta(D^\perp,j_1;n_2,s) = B_{j_1,0}(D^\perp)Q_0(0;n_2,s)n_2 - B_{j_1,0}(D^\perp)Q_1(0;n_2,s) - B_{j_1,1}(D^\perp)Q_1(1;n_2,s)$. Note that $Q_0(0;n_2,s) = 1$, $Q_1(0;n_2,s) = n_2(s-1)s^{-1}$ and $Q_1(1;n_2,s) = -s^{-1}$. Then

$$K_{t,1}(D) = \sum_{k=0}^{t}\binom{t}{k}(-1)^k n^{t-k}\sum_{j_1=0}^{n} Q_k(j_1;n,s)$$
$$\times [B_{j_1,0}(D^\perp)n_2 s^{-1} + B_{j_1,1}(D^\perp)s^{-1}]$$
$$= \sum_{j_1=0}^{n}\sum_{k=0}^{t}\binom{t}{k}(-1)^k n^{t-k} Q_k(j_1;n,s)$$
$$\times [B_{j_1,0}(D^\perp)n_2 s^{-1} + B_{j_1,1}(D^\perp)s^{-1}]$$
$$= \sum_{j_1=0}^{n} c_t(j_1;n,s)(s-1)^{-1}[B_{j_1,0}(D^\perp)n_2 + B_{j_1,1}(D^\perp)]s^{-1}.$$

Then (12) follows from (A.1) and the fact that $c_t(j_1;n,s) = 0$ when $j_1 > t$. □

PROOF OF LEMMA 3. By (8) and (15),

$$K_{t,1}(D) = N^{-2}\sum_{i=1}^{N}\sum_{j=1}^{N}[\delta_{ij}(D_T)]^t L_{p-1}$$
$$+ N^{-2}\sum_{m=1}^{s^p}\sum_{i\in B_m}\sum_{j\in B_m}[\delta_{ij}(D_T)]^t(L_p - L_{p-1})$$
$$= L_{p-1}K_{t,0}(D) + (L_p - L_{p-1})s^{-2p}\sum_{m=1}^{s^p} K_t(D_m).$$

Then the result follows from (16). □

PROOF OF LEMMA 4. Let $N_1 = s^{n-k-p}$ and $J_t(D_1) = \sum_{1 \le i < j \le N_1}[\delta_{ij}(D_1)]^t/[N_1(N_1-1)/2]$ for $t \ge 0$. It is easy to verify that, for $t \ge 0$,

$$K_t(D_1) = N_1^{-1}[(N_1-1)J_t(D_1) + n^t]. \tag{A.2}$$

(i) Recall that for an $(s^{n-k} : s^p)$ RME design $D$, $K_{1,0}(D) = s^{-1}n$ and $K_{1,1}(D) = s^{-2}nL_p$. Then, by Lemma 3, $K_1(D_1) = s[K_{1,1}(D) - L_{p-1}K_{1,0}(D)] = s^{-1}n$. By (A.2), $J_1(D_1) = (N_1-1)^{-1}[N_1 K_1(D_1) - n] = n(N_1-s)[(N_1-1)s]^{-1} = J$. On the other hand, Xu [22] showed that $J_1(D_1) \ge J$, with equality if and only if $D_1$ is an $OA(N_1, n, s, 1)$. Therefore, $D_1$ must be an $OA(N_1, n, s, 1)$.

(ii) Since the number of coincidences, $\delta_{ij}(D_1)$, must be an integer, it is easy to verify that, given $J_1(D_1) = J$, $J_2(D_1)$ achieves the minimum value of $J^2 + \eta(1-\eta)$ when all $\delta_{ij}(D_1)$, $i < j$, take on only one of the two values, $\lfloor J \rfloor$ and $\lfloor J \rfloor + 1$, where $\lfloor x \rfloor$ is the largest integer that does not exceed $x$. Then the result follows from (A.2).

(iii) By Lemma 3, $K_{2,1}(D) = L_{p-1}K_{2,0}(D) + s^{-1}K_2(D_1)$. The result follows from (ii) and the fact that $K_{2,0}(D) = s^{-2}n(n+s-1)$. $\square$

PROOF OF LEMMA 5. By the binomial theorem,

$$K_t(D_{\text{un}}) = N^{-2}\sum_{i=1}^{N}\sum_{j=1}^{N}[\delta_{ij}(D_T) + \delta_{ij}(D_F)]^t$$

$$= N^{-2}\sum_{i=1}^{N}\sum_{j=1}^{N}\sum_{r=0}^{t}\binom{t}{k}[\delta_{ij}(D_T)]^{t-r}[\delta_{ij}(D_F)]^r.$$

By (8), (15) and (16),

$$K_t(D_{\text{un}}) = N^{-2}\sum_{r=0}^{t}\binom{t}{k}\sum_{i=1}^{N}\sum_{j=1}^{N}[\delta_{ij}(D_T)]^{t-r}L_{p-1}^r$$

$$+ N^{-2}\sum_{r=0}^{t}\binom{t}{k}\sum_{m=1}^{s^p}\sum_{i \in B_m}\sum_{j \in B_m}[\delta_{ij}(D_T)]^{t-r}(L_p^r - L_{p-1}^r)$$

$$= \sum_{r=0}^{t}\binom{t}{k}L_{p-1}^r K_{t-r,0}(D) + \sum_{r=0}^{t}\binom{t}{k}(L_p^r - L_{p-1}^r)s^{-2p}\sum_{m=1}^{s^p}K_{t-r}(D_m)$$

$$= K_{t,0}(D) + \sum_{r=1}^{t}\binom{t}{k}L_{p-1}^r K_{t-r,0}(D)$$

$$+ s^{-p}\sum_{r=1}^{t}\binom{t}{k}(L_p^r - L_{p-1}^r)K_{t-r}(D_1).$$

Then the result follows from Lemma 3 with some algebra. □

PROOF OF LEMMA 6. First $A_{3,0}(D)$, $A_{2,1}(D)$ and $A_{4,0}(D)$ are minimized sequentially because $D$ has MA with respect to $W_{\text{scf}}$. Next, among designs with the same values of $A_{3,0}(D)$, $A_{2,1}(D)$ and $A_{4,0}(D)$, $A_{5,0}(D)$ is minimized because $D$ has MA with respect to $W_1$, and $A_{3,1}(D)$ is minimized because $D$ has MA with respect to $W_{\text{scf}}$. Continuing this type of argument shows that $D$ has MA with respect to $W_2$. □

PROOF OF LEMMA 7. First $3A_{3,0}(D) + A_{2,1}(D) = (3 - \alpha)A_{3,0}(D) + [\alpha A_{3,0}(D) + A_{2,1}(D)]$ is minimized because both $A_{3,0}(D)$ and $\alpha A_{3,0}(D) + A_{2,1}(D)$ are minimized. For designs with the same value of $3A_{3,0}(D) + A_{2,1}(D)$, they must have the same values of $A_{3,0}(D)$ and $A_{2,1}(D)$. Then $A_{4,0}(D)$ is minimized among designs with the minimum of $3A_{3,0}(D) + A_{2,1}(D)$ because $D$ has MA with respect to $W_{\text{scf}}$. Among designs with the same values of $3A_{3,0}(D) + A_{2,1}(D)$ and $A_{4,0}(D)$, $A_{5,0}(D)$ is minimized because $D$ has MA with respect to $W_2$ and $A_{3,1}(D)$ is minimized because $D$ has MA with respect to $W_{\text{scf}}$; therefore, $10A_{5,0}(D) + A_{3,1}(D)$ is also minimized. For designs with the same values of $3A_{3,0}(D) + A_{2,1}(D)$, $A_{4,0}(D)$ and $10A_{5,0}(D) + A_{3,1}(D)$, they must have the same $A_{5,0}(D)$ and $A_{3,1}(D)$ values. Continuing this type of argument shows that $D$ has MA with respect to $W_{\text{cc}}$. □

PROOF OF THEOREM 3. Note that $D$ achieves the lower bound in Theorem 2 if and only if $D_1$ achieves the lower bound in Lemma 4(ii). When the latter lower bound is achieved, $K_2(D_1)$ is minimized and $K_t(D_1)$ is uniquely determined for $t \geq 3$. By Lemma 3, $K_{t,1}(D)$ is determined by $K_{t,0}(D)$ for $t \geq 3$. Because $D_T$ has MA, by Lemma 2, $K_{3,0}(D), K_{4,0}(D), \ldots, K_{n,0}(D)$ are minimized sequentially. Then any combined sequence of $(K_{3,0}(D), K_{4,0}(D), \ldots, K_{n,0}(D))$ and $(K_{2,1}(D), K_{3,1}(D), \ldots, K_{n,1}(D))$ is also minimized sequentially as long as $K_{t,1}(D)$ is minimized after $K_{t,0}(D)$ for $t = 2, \ldots, n$. Hence, $D$ has minimum moment aberration and MA with respect to $W_{\text{scf}}$, $W_1$ and $W_2$. Finally, because $A_{2,1}(D)$ is minimized among all possible designs, by Lemma 7, $D$ has MA with respect to $W_{\text{cc}}$. □

PROOF OF THEOREM 4. *Necessity.* When $p = n - k - 1$, by Lemma 4(i), the principal block $D_1$ is an $OA(s, n, s, 1)$. Then it must contain a row of all zeros and other $s - 1$ rows without zeros.

*Sufficiency.* Let $u = (u_1, \ldots, u_n)$ be a row vector of $D_T$ without zeros. Because none of $u_i$ is zero, the linear equation $\sum_{i=1}^{n-k} x_i u_i = 0$ has $s^{n-k-1}$ solutions over $GF(s)$. Let $F$ be an $(n - k) \times L_{n-k-1}$ matrix, where the columns correspond to the solutions with the first nonzero element being unity. Clearly, $F$ has rank $n - k - 1$, and it is an $(n - k - 2)$-flat in $PG(n - k - 1, s)$. On the other hand, $D_T$ is an $[n, n - k]$ linear code over $GF(s)$. Let $T = (t_{ij})$ be the $(n-k) \times n$ generator matrix of $D_T$. We need to show that $T$ and $F$ have no columns in common so that the resulting blocked design $D = (D_T, D_F)$ is an RME design.

Without loss of generality, let $T = [I_{n-k}, E]$, where $I_{n-k}$ is the $n-k$ identity matrix and $E$ is an $(n-k) \times k$ matrix. Because the row vectors of $T$ form a basis for $D_T$, $u$ can be uniquely represented as a linear combination of the row vectors of $T$. Then it is clear that $\sum_{i=1}^{n-k} t_{ij} u_i = u_j \neq 0$ for $j = 1, \ldots, n$. This proves that $T$ and $F$ have no columns in common. $\square$

PROOF OF THEOREM 5. We only need to prove the MA optimality. Following the argument preceding the theorem, we can write an $(s^{n-k} : s^{n-k-1})$ RME design as $D = (D_T, D_F)$, where $D_T$ is a projection design of $\tilde{H}_{n-k}$. Recall that the principal block $D_1$ is an $OA(s, n, s, 1)$; therefore, each level appears exactly once in each column. It is evident that $\delta_{ij}(D_1) = 0$ when $i \neq j$ and $\delta_{ij}(D_1) = n$ when $i = j$; hence, $K_t(D_1) = s^{-1} n^t$ for $t > 0$. By Lemma 3, $K_{t,1}(D)$ is determined by $K_{t,0}(D)$ for $t > 0$. By (11) and (12), $A_{t,1}(D)$ is determined by $A_{1,0}(D), \ldots, A_{t,0}(D)$ uniquely. Thus, to sequentially minimize the sequences in (1), (2), (3) and (4), it is sufficient to sequentially minimize $A_{1,0}(D), \ldots, A_{n,0}(D)$. Then the result follows from the condition that $D_T$ has MA among all projection designs of $\tilde{H}_{n-k}$. $\square$

PROOF OF THEOREM 6. Given $K_{3,0}(D), K_{4,0}(D), \ldots, K_{t,0}(D)$, $3 \leq t \leq n$, by Lemma 5, sequentially minimizing $K_{2,1}(D), K_{3,1}(D), \ldots, K_{t-1,1}(D)$ is equivalent to sequentially minimizing $K_3(D_{\mathrm{un}}), K_4(D_{\mathrm{un}}), \ldots, K_t(D_{\mathrm{un}})$. Because $D_T$ has MA, by Lemma 2, $K_{3,0}(D), K_{4,0}(D), \ldots, K_{n,0}(D)$ are minimized sequentially. Because $D_{\mathrm{un}}$ has MA, $K_3(D_{\mathrm{un}}), K_4(D_{\mathrm{un}}), \ldots, K_n(D_{\mathrm{un}})$ are minimized sequentially. Then any combined sequence of $(K_{3,0}(D), K_{4,0}(D), \ldots, K_{n,0}(D))$ and $(K_{2,1}(D), K_{3,1}(D), \ldots, K_{n,1}(D))$ is also minimized sequentially as long as $K_{t-1,1}(D)$ is minimized after $K_{t,0}(D)$ for $t = 3, \ldots, n$. Hence, $D$ has minimum moment aberration and MA with respect to $W_{\mathrm{scf}}$, $W_1$ and $W_2$. By Lemma 5, $K_3(D_{\mathrm{un}}) = K_{3,0}(D) + 3K_{2,1}(D) + \mathit{constant}$; therefore, $K_{3,0}(D) + 3K_{2,1}(D)$ is minimized and, by (13) and (14), $A_{3,0}(D) + A_{2,1}(D)$ is minimized. Then, by Lemma 7, $D$ has MA with respect to $W_{\mathrm{cc}}$. $\square$

PROOF OF THEOREM 7. Because $D_T$ has MA, $A_{3,0}(D), A_{4,0}(D), \ldots, A_{n,0}(D)$ are minimized sequentially. Note that $A_{2,1}(D), A_{3,1}(D)$, and so on are minimized in $W_1$ or $W_2$ no sooner than in $W_{\mathrm{scf}}$. Therefore, if $D$ has MA with respect to $W_{\mathrm{scf}}$, it must have MA with respect to $W_1$ and $W_2$. The MA $W_{\mathrm{cc}}$ optimality follows from Lemma 7. $\square$

**Acknowledgments.** The author thanks a Co-Editor, an Associate Editor and two referees for their valuable comments.

Department of Statistics
University of California
Los Angeles, California 90095-1554
USA
E-mail: hqxu@stat.ucla.edu